\font\we=cmb10 at 14.4truept
\font\li=cmb10 at 12truept
~\vskip 1.0cm
\centerline {\we Constructions of Non-Abelian Zeta Functions for Curves}
\vskip 1.0cm
\centerline {\li Lin WENG}
\centerline {\bf Graduate School of Mathematics, Nagoya University, Japan}
\vskip 0.45cm
In this paper, we initiate a geometrically oriented study of local and 
global non-abelian zeta
functions for curves. This consists of two parts: construction and 
justification.

For the construction,  we first use moduli spaces of semi-stable
bundles to introduce a new type of zeta functions for curves defined over 
finite fields. Then,
we prove that these new zeta
functions are indeed rational and satisfy the functional equation, based
on vanishing theorem, duality, Riemann-Roch theorem for cohomology of 
semi-stable vector bundles.
With all this,   we next introduce certain global non-abelian zeta 
functions for curves
defined over number fields, via the  Euler product formalism. Finally we 
establish the convergence
of our Euler products, from the Clifford Lemma, an ugly yet explicit 
formula for local non-abelian
zeta functions,  a  result of (Harder-Narasimhan) Siegel about quadratic
forms, and Weil's theorem on Riemann Hypothesis of Artin zeta functions.

As for justification, surely, we check that when only line bundles are
involved, (so  moduli spaces of semi-stable bundles are nothing but the 
standard  Picard groups),
our (new) zeta functions coincide with the classical  Artin zeta functions 
for curves over finite
fields  and Hasse-Weil zeta functions for curves over number fields 
respectively. Moreover we compute
the rank two zeta functions for genus two curves by studying the so-called 
non-abelian Brill-Noether
loci and their infinitesimal structures. This is indeed a  qutie 
interesting, and in general, should
be a very  important aspect of the theory:  We not only need to precisely 
describe all of the
Brill-Noether loci but the so-called associated infinitesmal structures 
attached to all Seshadri
equivalence classes, in which Weierstrass points appear naturally.

This work is motivated by our studies on new non-abelian zeta functions for 
number fields and
Tamagawa measures associated to the Weil-Petersson metrics on moduli spaces 
of stable vector
bundles. So related conjectures, or better, working hypothesis, are 
proposed. We hope that our
non-abelian zeta functions, which do not really follow the present style of 
the theory of zeta
functions, are acceptable and hence play a certain role in exploring the 
non-abelian aspect of
arithmetic of curves.

\vskip 0.45cm
\centerline {\li Chapter 1. Local Non-Abelian Zeta Functions for Curves}
\vskip 0.30cm
In this chapter, we introduce our non-abelian
zeta functions for curves defined over number fields. Basic properties such 
as meromorphic
extensions, rationality and functional equations are established.
\vskip 0.30cm
\centerline {\bf 1.1. Moduli Spaces of Semi-Stable Bundles}
\vskip 0.30cm
\noindent
{\bf 1.1.1. Semi-Stable Bundles.} Let $C$ be a regular, reduced and 
irreducible projective curve
defined over an algebraically closed field
$\bar k$. Then according to Mumford [M], a vector bundle $V$ on $C$ is 
called semi-stable (resp.
stable) if for any proper subbundle $V'$ of $V$,  $$\mu(V'):={{d(V')}\over 
{r(V')}}\leq {\rm
(resp.}\leq {\rm )}{{d(V)}\over {r(V)}}=:\mu(V).$$ Here $d$ denotes the 
degree and $r$ denotes the
rank.

\noindent
{\bf Proposition.} {\it Let $V$ be a vector bundle over $C$. Then

\noindent
(a) ([HN]) there exists a unique
filtration of subbundles of $V$, the  Harder-Narasimhan filtration of $V$,
$$\{0\}=V_0\subset V_1\subset V_2\subset\dots\subset V_{s-1}\subset V_s=V$$
such that for $1\leq i\leq s-1$, $V_i/V_{i-1}$ is semi-stable and
$\mu(V_i/V_{i-1}>\mu(V_{i+1}/V_i);$

\noindent
(b) (see e.g. [Se]) if moreover $V$ is semi-stable,
there exists a
filtration of subbundles of $V$, a Jordan-H\"older filtration of $V$,
$$\{0\}=V^{t+1}\subset V^t\subset\dots\subset V^{1}\subset V^0=V$$
such that for all $0\leq i\leq t$, $V^i/V^{i+1}$ is stable and 
$\mu(V^i/V^{i+1})=\mu (V)$. Moreover,
the associated graded bundle ${\rm Gr}(V):=\oplus_{i=0}^tV^i/V^{i+1}$, the
(Jordan-H\"older) graded bundle of
$V$, is determined uniquely by $V$.}
\vskip 0.30cm
\noindent
{\bf 1.1.2. Moduli Space of Stable Bundles.} Following Seshadri, two 
semi-stable vector bundles $V$
and $W$ are called $S$-equivalent, if their associated Jordan-H\"older 
graded bundles are
isomorphic, i.e.,
${\rm Gr}(V)\simeq {\rm Gr}(W)$. Applying Mumford's general result on 
geometric invariant theory,
Seshadri proves the following

\noindent
{\bf Theorem.} ([Se]) {\it Let $C$ be a regular, reduced, irreducible 
projective curve of genus
$g\geq 2$ defined over an algebraically closed field. Then over the set 
${\cal M}_{C,r}(d)$ (resp.
${\cal M}_{C,r}(L)$) of $S$-equivalence classes of rank $r$ and degree $d$ 
(resp. rank $r$ and
determinant $L$) semi-stable vector bundles over $C$, there is a natural 
normal, projective
$(r^2(g-1)+1)$-dimensional (resp. $(r^2-1)(g-1)$-dimensional) algebraic 
variety structure.}

\noindent
{\it Remark.} In this paper, we always assume that the genus of $g$ is at 
least 2. For  elliptic
curves, whose associated moduli spaces are very special, please see [We2].
\vskip 0.30cm
\noindent
{\bf 1.1.3. Rational Points.} Now  assume that $C$ is defined over a finite 
field $k$. Naturally
we may talk about
$k$-rational bundles over $C$, i.e., bundles which are defined over $k$. 
Moreover, from geometric
invariant theory, projective varieties
${\cal M}_{C,r}(d)$ are defined over a certain finite extension of $k$; and 
if $L$ itself is
defined over $k$, the same holds for ${\cal M}_{C,r}(L)$. Thus it makes 
sense to talk about
$k$-rational points of these moduli spaces too. The relation between these 
two types of rationality
is given by  Harder and Narasimhan based on a discussion about Brauer groups:

\noindent
{\bf Proposition.} ([HN]) {\it Let $C$ be a regular, reduced, irreducible 
projective curve of genus
$g\geq 2$ defined over a finite field $k$. Then there exists a finite field 
${\bf F}_q$ such that
for all
$d$ (resp. all $k$-rational line bundles $L$), the subset of ${\bf 
F}_q$-rational points of
${\cal M}_{C,r}(d)$ (resp. ${\cal M}_{C,r}(L)$) consists exactly of all 
$S$-equivalence classes of
${\bf F}_q$-rational bundles in ${\cal M}_{C,r}(d)$ (resp. ${\cal
M}_{C,r}(L)$).}

 From now on, without loss of generality, we always assume that  finite fields
${\bf F}_q$ (with $q$ elements) satisfy the property in the Proposition. 
Also for simplicity, we
write ${\cal M}_{C,r}(d)$ (resp.
${\cal M}_{C,r}(L)$) for ${\cal M}_{C,r}(d)({\bf F}_q)$ (resp. ${\cal 
M}_{C,r}(L)({\bf
F}_q)$), the subset of ${\bf F}_q$-rational points, and call them moduli 
spaces by an abuse of
notations.
\vskip 0.45cm
\centerline {\bf 1.2. Local Non-Abelian Zeta Functions}
\vskip 0.30cm
\noindent
{\bf 1.2.1. Definition.} Let $C$ be a  regular, reduced, irreducible 
projective curve of genus $g\geq
2$ defined over the finite field ${\bf F}_q$ with $q$ elements. Define the 
{\it rank $r$
non-abelian zeta function} $\zeta_{C,r,{\bf F}_q}(s)$ by setting
$$\zeta_{C,r,{\bf F}_q}(s):=\sum_{V\in [V]\in {\cal M}_{C,r}(d), d\geq 
0}{{q^{h^0(C,V)}-1}\over
{\#{\rm Aut}(V)}}\cdot (q^{-s})^{d(V)},\qquad {\rm Re}(s)>1.$$

\noindent
{\bf Proposition.} {\it With the same notation as above, $\zeta_{C,1,{\bf 
F}_q}(s)$ is nothing but
the classical Artin zeta function for curve $C$. That is to say,
$$\zeta_{C,1,{\bf F}_q}(s)=\sum_{D\geq 0}{1\over {N(D)^{s}}}\qquad {\rm 
Re}(s)>1.$$ Here $D$ runs
over all effective divisors of $C$.}

\noindent
{\it Proof.} By definition, the classical Artin zeta function ([A]) for $C$ 
is given by
$$\zeta_C(s):=\sum_{D\geq 0}{1\over {N(D)^{s}}}.$$ Here $N(D)=q^{d(D)}$ with
$d(\Sigma_Pn_PP)=\Sigma_Pn_Pd(P)$. Thus by first grouping effective 
divisors according to their
rational equivalence classes ${\cal D}$, then taking the sum on effective 
divisors in the same class,
we obtain $$\zeta_C(s)=\sum_{\cal D}\sum_{D\in {\cal D},D \geq 0}{1\over 
{N(D)^{s}}}.$$
Clearly, $$\sum_{D\in {\cal D}, D \geq 0}{1\over 
{N(D)^{s}}}={{q^{h^0(C,{\cal D})}-1}\over
{q-1}}\cdot (q^{-s})^{d({\cal D})}.$$ Therefore, $$\zeta_C(s)=\sum_{L\in {\rm
Pic}^d(C), d\geq 0}{{q^{h^0(C,L)}-1}\over {\#{\rm Aut}(L)}}\cdot 
(q^{-s})^{d(L)}$$ due to the
fact that
${\rm Aut}(V)\simeq {\bf F}_q^*$.
\vskip 0.30cm
\noindent
{\bf 1.2.2. Convergence and Rationality.} Clearly, we must show that for 
general $r$, the infinite
summation in the definition of our non-abelian zeta function 
$\zeta_{C,r,{\bf F}_q}(s)$ converges
when ${\rm Re}(s)>1$. For this, let us start with  the following vanishing 
result for semi-stable
bundles.

\noindent
{\bf Lemma 1.} {\it Let $V$ be a degree $d$ and rank $r$ semi-stable vector 
bundle on $C$. Then

\noindent
(a) if $d\geq r(2g-2)+1$, $h^1(C,V)=0$;

\noindent
(b) if $d<0$, $h^0(C,V)=0$.}

\noindent
{\it Proof.} This is a direct consequence of the fact that if $V$ and $W$ 
are semi-stable bundles
with $\mu(V)>\mu(W)$, then $H^0(C,{\rm Hom}(V,W))=\{0\}$.

Thus, from definition,
$$\eqalign{\zeta_{C,r,{\bf F}_q}(s)=&
\sum_{V\in [V]\in {\cal M}_{C,r}(d), 0\leq d\leq r(2g-2)
}{{q^{h^0(C,V)}-1}\over {\#{\rm Aut}(V)}}\cdot (q^{-s})^{d(V)}\cr
&\qquad+\sum_{V\in [V]\in {\cal
M}_{C,r}(d), d\geq r(2g-2)+1}{{q^{d(V)-r(g-1)}-1}\over {\#{\rm 
Aut}(V)}}\cdot (q^{-s})^{d(V)}.\cr}$$
Clearly only finitely many terms appear in the first summation, so it 
suffices to show that when
${\rm Re}(s)>1$, the second term converges. For this purpose, we introduce 
the so-called
Harder-Narasimhan numbers
$$\beta_{C,r,{\bf F}_q}(d):=\sum_{V\in [V]\in {\cal M}_{C,r}(d)}{1\over 
{\#{\rm Aut}(V)}}.$$

\noindent
{\bf Lemma 2.} {\it With the same notation as above, for all $n\in {\bf Z}$,
$$\beta_{C,r,{\bf F}_q}(d+rn)=\beta_{C,r,{\bf F}_q}(d).$$}

\noindent
{\it Proof.} This comes from the following two facts: (1) there is a degree 
one ${\bf F}_q$-rational
line bundle $A$ on $C$; and (2) ${\rm Aut}(V)\simeq
{\rm Aut}(V\otimes A^{\otimes n})$ and $d(V\otimes A^{\otimes n})=d(V)+rn$.

Therefore,  the second summation becomes
$$\eqalign{~&\sum_{i=1}^{r}\beta_{C,r,{\bf
F}_q}(i)\sum_{n=2g-2}^\infty\Big(q^{nr+i-r(g-1)}-1\Big)\cdot
(q^{-s})^{nr+i}\cr
=&\sum_{i=1}^r\beta_{C,r,{\bf F}_q}(i) \cdot (q^{-s})^i\cdot
\Big(q^{i-r(g-1)}\cdot {{q^{(1-s)\cdot r(2g-2)}}\over {1-q^{(1-s)\cdot r}}}
-{{q^{(-s)\cdot r(2g-2)}}\over {1-q^{(-s)\cdot r}}}\Big),\cr}$$
provided that $|q^{-s}|<1$. Thus we have proved the following

\noindent
{\bf Proposition.} {\it The non-abelian zeta function $\zeta_{C,r,{\bf 
F}_q}(s)$ is well-defined
for ${\rm Re}(s)>1$, and admits a meromorphic extension to the whole 
complex $s$-plane.}

Moreover, if we  set
$t:=q^{-s}$ and introduce the non-abelian $Z$-function of $C$ by setting
$$\zeta_{C,r,{\bf F}_q}(s)=:Z_{C,r,{\bf F}_q}(t):=\sum_{V\in [V]\in {\cal 
M}_{C,r}(d),d\geq
0}{{q^{h^0(C,V)}-1}\over {\#{\rm Aut}(V)}}\cdot t^{d(V)}, \qquad |t|<1.$$
Then the above calculation implies that
$$Z_{C,r,{\bf F}_q}(t)
=\sum_{d=0}^{r(2g-2)}\Big(\sum_{V\in [V]\in {\cal
M}_{C,r}(d)}{{q^{h^0(C,V)}-1}\over {\#{\rm Aut}(V)}}\Big)\cdot t^{d}
+\sum_{i=1}^r\beta_{C,r,{\bf
F}_q}(i)\cdot \Big({{q^{r(g-1)+i}}\over {1-q^rt^r}}-{{1}\over 
{1-t^r}}\Big)\cdot t^{r(2g-2)+i}.$$
Therefore, there exists a polynomial $P_{C,r,{\bf F}_q}(s)\in {\bf Q}[t]$ 
such that
$$Z_{C,r,{\bf F}_q}(t)={{P_{C,r,{\bf F}_q}(t)}\over {(1-t^r)(1-q^rt^r)}}.$$ 
In this way, we have
established the following

\noindent
{\bf Rationality.} {\it Let $C$ be a regular, reduced irreducible 
projective curve defined over ${\bf
F}_q$ with $Z_{C,r,{\bf F}_q}(t)$ the rank $r$ non-abelian $Z$-function. 
Then, there exists a
polynomial
$P_{C,r,{\bf F}_q}(s)\in {\bf Q}[t]$ such that
$$Z_{C,r,{\bf F}_q}(t)={{P_{C,r,{\bf F}_q}(t)}\over {(1-t^r)(1-q^rt^r)}}.$$}

\noindent
{\bf 1.2.3. Functional Equation.} Besides the fact that
$P_{C,r,{\bf F}_q}(t)\in {\bf Q}[t]$, we know  very little about this 
polynomial. To understand
$P_{C,r,{\bf F}_q}(s)$ better, as well as for the theoretical purpose, we 
next study the functional
equation for  rank $r$ zeta functions. As a preparation, here we introduce 
the rank $r$
non-abelian $\xi$-function
$\xi_{C,r,{\bf F}_q}(s)$ by setting
$$\xi_{C,r,{\bf F}_q}(s):=\zeta_{C,r,{\bf F}_q}(s)\cdot (q^{s})^{r(g-1)}.$$ So
$$\xi_{C,r,{\bf F}_q}(s)=\sum_{V\in [V]\in {\cal M}_{C,r}(d),d\geq 
0}{{q^{h^0(C,V)}-1}\over {\#{\rm
Aut}(V)}}\cdot (q^{-s})^{\chi(C,V)}, \qquad {\rm Re}(s)>1,$$
where $\chi(C,V)$ denotes the Euler-Poincar\'e characteristic of $V$.

\noindent
{\bf Functional Equation.} {\it Let $C$ be a regular, reduced irreducible 
projective curve defined
over
${\bf F}_q$ with $\xi_{C,r,{\bf F}_q}(s)$ its associated rank $r$ 
non-abelian  $\xi$-function. Then,
$$\xi_{C,r,{\bf F}_q}(s)=\xi_{C,r,{\bf F}_q}(1-s).$$}

Before proving  the functional equation, we give the following

\noindent
{\bf Corollary.} {\it With the same notation as above,

\noindent
(a) $P_{C,r,{\bf F}_q}(t)\in {\bf Q}[t]$ is a degree $2rg$ polynomial;

\noindent
(b) Denote all reciprocal roots of $P_{C,r,{\bf F}_q}(t)$ by 
$\omega_{C,r,{\bf F}_q}(i),
i=1,\dots, 2rg$. Then after a suitable rearrangement,
$$\omega_{C,r,{\bf F}_q}(i)\cdot \omega_{C,r,{\bf F}_q}(2rg-i)=q,\qquad 
i=1,\dots,rg;$$

\noindent
(c) For each $m\in {\bf Z}_{\geq 1}$, there exists a rational number 
$N_{C,r,{\bf F}_q}(m)$ such
that
$$Z_{r,C,{\bf F}_q}(t)=P_{C,r,{\bf F}_q}(0)\cdot\exp\Big(\sum_{m=1}^\infty 
N_{C,r,{\bf
F}_q}(m){{t^m}\over m}\Big).$$ Moreover, $$N_{C,r,{\bf
F}_q}(m)=\cases{r(1+q^m)-\sum_{i=1}^{2rg}\omega_{C,r,{\bf F}_q}(i)^m,& if 
$r\ |m$;\cr
-\sum_{i=1}^{2rg}\omega_{C,r,{\bf F}_q}(i)^m,& if
$r\not| m$;\cr}$$

\noindent
(d) For any  $a\in {\bf Z}_{>0}$, denote by $\zeta_{a}$ a primitive $a$-th 
root of unity and set
$T=t^a$.  Then
$$\prod_{i=1}^aZ_{C,r}(\zeta_{a}^it)=(P_{C,r,{\bf F}_q}(0))^a\cdot
\exp\Big(\sum_{m=1}^\infty N_{r,C,{\bf F}_q}(ma){{T^m}\over m}\Big).$$}

\noindent
{\it Proof.} (a) and (b) are direct consequences of the functional 
equation, while (c) and (d) are
direct consequences of (a), (b) and the following well-known relations
$$\sum_{i=1}^a(\zeta_{a}^i)^m=\cases{a,& if $a\ |m$,\cr
0, & if $a\not|m$.\cr}$$

\noindent
{\bf 1.2.4. Proof of Functional Equation.} To understand the structure of 
the functional equation
clearly, we need to decompose the non-abelian $\xi$-function for curves.
For this purpose, first recall that the canonical line bundle $K_C$ of $C$ 
is defined over ${\bf
F}_q$. Thus, for all
$n\in {\bf Z}$, we obtain the following natural ${\bf F}_q$-rational 
isomorphisms:
$$\matrix{ {\cal M}_r(L)&\to& {\cal M}_r(L\otimes K_C^{\otimes nr});&\qquad&
{\cal M}_r(L)&\to& {\cal M}_r(L^{\otimes -1}\otimes K_C^{\otimes nr})\cr
[V]&\mapsto&[V\otimes K_C^{\otimes n}];&\qquad&[V]&\mapsto& [V^\vee\otimes 
K_C^{\otimes
n}],\cr}$$  where  $V^\vee$ denotes the dual of $V$. Next, introduce the union
$${\cal M}_{C,r}^L:=\cup_{n\in {\bf Z}}\Big({\cal M}_r(L\otimes
K_C^{\otimes nr})\cup {\cal M}_r(L^{\otimes -1}\otimes K_C^{\otimes
nr})\Big).$$ With this, clearly, we may and indeed always assume that
$$0\leq d(L)\leq r(g-1).$$

Further, introduce the partial non-abelian zeta function $\xi_{C,r,{\bf 
F}_q}^L(s)$ by setting
$$\xi_{C,r,{\bf F}_q}^L(s):=\sum_{V\in [V]\in {\cal 
M}_{C,r}^L}{{q^{h^0(C,V)}-1}\over {\#{\rm
Aut}(V)}} \cdot\big(q^{-s}\big)^{\chi(C,V)},\qquad {\rm Re}(s)>1.$$ 
Clearly, then
$$\xi_{C,r,{\bf F}_q}(s)=\sum_L\xi_{C,r,{\bf F}_q}^L(s)$$ where $L$ runs 
over all line bundles
appeared in the following (disjoint) union
$$\cup_{d\in {\bf Z}}{\cal M}_{C,r}(d)=\cup_{L,0\leq d(L)\leq r(g-1)}{\cal 
M}_{C,r}^L.$$

\noindent
{\it Remark.} Here we reminder the reader that the vanishing result of 
Lemma 1.2.2.1 has been used.
 
Therefore, to prove the functional equation for $\xi_{C,r,{\bf F}_q}(s)$, 
it suffices to show
$$\xi_{C,r,{\bf F}_q}^L(s)=\xi_{C,r,{\bf F}_q}^L(1-s).$$
For this, we have the following

\noindent
{\bf Theorem.} {\it  For ${\rm Re}(s)>1$,
$$\eqalign{~&\xi_{C,r,{\bf F}_q}^L(s)\cr
=&{1\over
2}\sum_{V\in [V]\in {\cal M}_{C,r}^L;0\leq d(V)\leq
r(2g-2)}{{q^{h^0(C,V)}}\over {\#{\rm
Aut}(V)}}\cdot\Big[(q^{-s})^{\chi(C,V)}+
(q^{s-1})^{\chi(C,V)}\Big]\cr
&+\Big[{{q^{(1-s)\cdot(d(L)-r(g-1))}}\over
{q^{(s-1)\cdot r(2g-2)}-1}}
+{{q^{s\cdot(d(L)-r(g-1))}}\over
{q^{(-s)\cdot r(2g-2)}-1}}+{{q^{(s-1)\cdot(d(L)-r(g-1))}}\over
{q^{(s-1)\cdot r(2g-2)}-1}}
+{{q^{(-s)\cdot(d(L)-r(g-1))}}\over
{q^{(-s)\cdot r(2g-2)}-1}}\Big]\cdot \beta_{C,r,{\bf 
F}_q}(L).\cr}\eqno(*)$$ Here
$\beta_{C,r,{\bf F}_q}(L):=\sum_{V\in [V]\in {\cal M}_{C,r}(L)}{1\over{\#{\rm
Aut}(E)}}$ denotes the Harder-Narasimhan number. In particular,

\noindent
(a) $\xi_{C,r,{\bf F}_q}^L(s)$ satisfies the functional equation
$\qquad\xi_{C,r,{\bf F}_q}^L(s)=\xi_{C,r,{\bf F}_q}^L(1-s);$

\noindent
(b) the Harder-Narasimhan number $\beta_{C,r,{\bf F}_q}(L)$ is given by the 
leading
term of the  singularities of $\xi_{C,r,{\bf F}_q}^L(s)$  at $s=0$ and $s=1$.}

\noindent
{\it Remark.} In this way, we use non-abelian zeta functions to evaluate 
the so-called
Narasimhan-Harder numbers. Thus, the Betti numbers for moduli spaces of 
stable bundles may also be
read from our non-abelian zeta functions. The global version of this will 
be discussed in 2.2.3
below.

\noindent
{\it Proof.} It suffices to prove (*). For this, set
$$I(s)=\sum_{V\in [V]\in {\cal M}_{C,r}^L;0\leq d(V)\leq
r(2g-2)}{{q^{h^0(C,V)}}\over {\#{\rm
Aut}(V)}}\cdot (q^{-s})^{\chi(C,V)}$$ and
$$II(s)=\sum_{V\in [V]\in {\cal M}_{C,r}^L;=d(V)>
r(2g-2)}{{q^{h^0(C,V)}}\over {\#{\rm
Aut}(V)}}\cdot (q^{-s})^{\chi(C,V)}-\sum_{V\in [V]\in {\cal 
M}_{C,r}^L;d(V)\geq 0}{{1}\over {\#{\rm
Aut}(V)}}\cdot (q^{-s})^{\chi(C,V)}.$$ Thus,
$$\xi_{C,r,{\bf F}_q}^L(s)=I(s)+II(s).$$  So it suffices to show the following

\noindent
{\bf Lemma.} {\it With the same notation as above,

\noindent
(a) $I(s)={1\over
2}\sum_{V\in [V]\in {\cal M}_{C,r}^L;0\leq d(V)\leq
r(2g-2)}{{q^{h^0(C,V)}}\over {\#{\rm
Aut}(V)}}\cdot\Big[(q^{-s})^{\chi(C,V)}+
(q^{s-1})^{\chi(C,V)}\Big];$ and

\noindent
(b) $$II(s)=\Big[{{q^{(1-s)\cdot(d(L)-r(g-1))}}\over
{q^{(s-1)\cdot r(2g-2)}-1}}
+{{q^{s\cdot(d(L)-r(g-1))}}\over
{q^{(-s)\cdot r(2g-2)}-1}}+{{q^{(s-1)\cdot(d(L)-r(g-1))}}\over
{q^{(s-1)\cdot r(2g-2)}-1}}
+{{q^{(-s)\cdot(d(L)-r(g-1))}}\over
{q^{(-s)\cdot r(2g-2)}-1}}\Big]\cdot \beta_{C,r,{\bf F}_q}(L).$$}

\noindent
{\it Proof.} (a) comes from  Riemann-Roch theorem and  Serre duality. Indeed,
$$\eqalign{I(s)
=&{1\over
2}\Big(\sum_{V\in [V]\in {\cal M}_{C,r}^L;0\leq d(E)\leq
r(2g-2)}{{q^{h^0(C,V)}}\over {\#{\rm
Aut}(V)}}\cdot(q^{-s})^{\chi(C,V)}\cr
&\qquad+\sum_{V^\vee\otimes
K_C\in {\cal M}_{C,r}^L;0\leq d(V^\vee\otimes K_C)\leq
r(2g-2)}{{q^{h^0(C,V^\vee\otimes K_C)}}\over {\#{\rm
Aut}(V^\vee\otimes K_C)}}\cdot(q^{-s})^{\chi(C,V^\vee\otimes
K_C)}\Big)\cr
=&{1\over
2}\sum_{V\in [V]\in {\cal M}_{C,r}^L;0\leq d(E)\leq
r(2g-2)}\Big[{{q^{h^0(C,V)}}\over {\#{\rm
Aut}(V)}}\cdot(q^{-s})^{\chi(C,V)}+{{q^{h^1(C,V^\vee\otimes K_C)}}\over {\#{\rm
Aut}(V^\vee\otimes K_C)}}\cdot(q^{1-s})^{\chi(C,V^\vee\otimes
K_C)}\Big]\cr
=&{1\over
2}\sum_{V\in [V]\in {\cal M}_{C,r}^L;0\leq d(V)\leq
r(2g-2)}{{q^{h^0(C,V)}}\over {\#{\rm
Aut}(V)}}\cdot\Big[(q^{-s})^{\chi(C,V)}+
(q^{s-1})^{\chi(C,V)}\Big].\cr}$$

As for (b), clearly by the vanishing result,
$$\eqalign{T_{r,L}(s)=&\sum_{V\in [V]\in {\cal M}_{C,r}^L;d(E)>
r(2g-2)}{1\over {\#{\rm
Aut}(E)}}\cdot(q^{1-s})^{\chi(C,V)}-\sum_{V\in [V]\in {\cal 
M}_{C,r}^L;d(E)\geq 0}
{1\over {\#{\rm Aut}(E)}}\cdot(q^{-s})^{\chi(C,V)}\cr
=&\Big(\sum_{V\in [V]\in {\cal M}_{C,r}(L\otimes K_C^{\otimes rn});
d(L)+rn(2g-2)> r(2g-2)}
{1\over {\#{\rm Aut}(V)}}\cdot(q^{1-s})^{\chi(C,E)}\cr
&\qquad-\sum_{V\in [V]\in {\cal M}_{C,r}(L^{-1}\otimes K_C^{\otimes rn});
-d(L)+rn(2g-2)\geq 0}
{1\over {\#{\rm Aut}(V)}}\cdot(q^{-s})^{\chi(C,E)}\Big)\cr
&+\Big(\sum_{V\in [V]\in {\cal M}_{C,r}(L^{-1}\otimes K_C^{\otimes rn});
-d(L)+rn(2g-2)> r(2g-2)}
{1\over {\#{\rm Aut}(V)}}\cdot(q^{1-s})^{\chi(C,E)}\cr
&\qquad-\sum_{V\in [V]\in {\cal M}_{C,r}(L\otimes K_C^{\otimes rn});
d(L)+rn(2g-2)>0}
{1\over {\#{\rm Aut}(V)}}\cdot(q^{-s})^{\chi(C,E)}\Big).\cr}$$
But $\chi(C,V)$ depends only on $d(V)$. Thus, accordingly,
$$\eqalign{II(s)
=&\Big[\Big(\sum_{n=1}^\infty (q^{1-s})^{d(L)+nr(2g-2)-r(g-1)}
-\sum_{n=1}^\infty(q^{-s})^{-d(L)+nr(2g-2)-r(g-1)}\Big)\cr
&+\Big(\sum_{n=2}^\infty(q^{1-s})^{-d(L)+nr(2g-2)-r(g-1)}
-\sum_{n=0}^\infty(q^{-s})^{d(L)+nr(2g-2)-r(g-1)}\Big)\Big]\cdot
\beta_{C,r}(L)\cr
=&\Big[{{q^{(1-s)\cdot(d(L)-r(g-1))}}\over
{q^{(s-1)\cdot r(2g-2)}-1}}
+{{q^{s\cdot(d(L)-r(g-1))}}\over
{q^{(-s)\cdot r(2g-2)}-1}}+{{q^{(s-1)\cdot(d(L)-r(g-1))}}\over
{q^{(s-1)\cdot r(2g-2)}-1}}
+{{q^{(-s)\cdot(d(L)-r(g-1))}}\over
{q^{(-s)\cdot r(2g-2)}-1}}\Big]\cdot \beta_{C,r,{\bf F}_q}(L).\cr}$$
This completes the proof of the lemma, and hence the Theorem and the
Functional Equation for rank $r$ zeta functions of curves.
\vskip 0.45cm
\centerline {\li Chapter 2. Global Non-Abelian Zeta Functions for Curves}
\vskip 0.30cm
In this chapter, we introduce new non-abelian zeta functions for curves 
defined over number fields
via the Euler product formalism. Thus the global construction is based on 
our study for the
non-abelian zeta of curves defined over finite fields in Chapter 1. Main 
result of this chapter is
about a convergence region of such Euler products. Key ingredients of the 
proof are a result of
(Harder-Narasimhan) Siegel, an ugly yet very precise formula for our local 
zeta functions,
Clifford Lemma for semi-stable vector bundles, and Weil's theorem on 
Riemann Hypothesis for
Artin zeta functions.
\vskip 0.30cm
\centerline {\bf 2.1. Preparations}
\vskip 0.30cm
\noindent
{\bf 2.1.1. Invariants $\alpha,\beta$ and $\gamma$.} Let $C$ be a regular, 
reduced, irreducible
projective curve of genus $g$ defined over the finite field
${\bf F}_q$ with $q$  elements. As in Chapter 1, we then get (the subset of 
${\bf F}_q$-rational
points of) the associated moduli spaces
${\cal M}_{E,r}({\cal L})$ and ${\cal M}_{C,r}(d)$. Recall that in Chapter 1,
motivated by a work of Harder-Narasimhan [HN], we, following Desale-Ramanan 
[DR],
define the Harder-Narasimhan numbers $\beta_{C,r,{\bf F}_q}(L), 
\beta_{C,r,{\bf F}_q}(d)$, which
are very useful in the discussion of our zeta functions.  Now  we introduce 
new invariants
for $C$ by setting
$$\alpha_{C,r,{\bf F}_q}(d):=\sum_{V\in [V]\in {\cal M}_{C,r}(d)({\bf
F}_q)}{{q^{h^0(C,V)}}\over {\#{\rm Aut}(V)}},\qquad\qquad
\gamma_{C,r,{\bf F}_q}(d):=\sum_{V\in [V]\in {\cal M}_{C,r}(d)({\bf
F}_q)}{{q^{h^0(C,V)}-1}\over {\#{\rm Aut}(V)}},$$
and similarly define $\alpha_{C,r,{\bf F}_q}(L)$ and $\gamma_{C,r,{\bf 
F}_q}(L)$.

\noindent
{\bf Lemma.} {\it With the same notation as above,

\noindent
(1) for $\alpha_{C,r,{\bf F}_q}(d)$,  $$\alpha_{C,r,{\bf 
F}_q}(d)=\cases{\beta_{C,r,{\bf
F}_q}(d);&if
$d< 0;$\cr
\alpha_{C,r,{\bf F}_q}(r(2g-2)-d)\cdot q^{d-r(g-1)},& if  $0\leq d\leq 
r(2g-2);$\cr
\beta_{C,r,{\bf F}_q}(d)\cdot q^{d-r(g-1)},&if
$d>r(2g-2);$\cr}$$

\noindent
(ii) for $\beta_{C,r,{\bf F}_q}(d)$,
$$\beta_{C,r,{\bf F}_q}(\pm d+rn)=\beta_{C,r,{\bf F}_q}(d)\qquad n\in {\bf 
Z};$$

\noindent
(iii) for $\gamma_{C,r,{\bf F}_q}(d)$, $$\gamma_{C,r,{\bf
F}_q}(d)=\alpha_{C,r,{\bf F}_q}(d)-\beta_{C,r,{\bf F}_q}(d).$$}

\noindent
{\it Proof.}  (iii) is simply the definition, while (ii) is a direct 
consequence of Lemma 1.2.2.2 and
the fact that ${\rm Aut}(V)\simeq {\rm Aut}(V^\vee)$ for a vector bundle 
$V$. So it suffices to
prove (i).

When $d<0$, the relation is deduced from the fact that $h^0(C,V)=0$ if
$V$ is a semi-stable vector bundle with strictly negative degree; when 
$0\leq d\leq r(2g-2)$, the
result comes from the Riemann-Roch and Serre duality; finally when 
$d>r(2g-2)$, the result
is a direct consequence of the Riemann-Roch and the fact that  $h^1(C,V)=0$ 
if $V$ is
a semi-stable vector bundle with degree strictly bigger  than $r(2g-2)$.

We here reminder the reader that this Lemma and  Lemma 1.2.2.2 tell us that all
$\alpha_{C,r,{\bf F}_q}(d),\beta_{C,r,{\bf F}_q}(d)$ and $\gamma_{C,r,{\bf 
F}_q}(d)$'s may be
calculated from $\alpha_{C,r,{\bf F}_q}(i),\beta_{C,r,{\bf F}_q}(j)$ with 
$i=0,\dots, r(g-1)$ and
$j=0,\dots,r-1$.

\noindent
{\bf 2.1.2. Asymptotic Behaviors of $\alpha,\beta$ and $\gamma$}

For later use, we here discuss the asymptotic behavior of $\alpha_{C,r,{\bf 
F}_q}(d)$,
$\beta_{C,r,{\bf F}_q}(d)$, and $\gamma_{C,r, {\bf F}_q}(0)$ when $q\to 
\infty$.
\vskip 0.30cm
\noindent
{\bf Proposition.} {\it With the same notation as above, when $q\to\infty$,

\noindent
(a) For all $d$, $$\beta_{C,r,{\bf
F}_q}(d)=O\Big(q^{r^2(g-1)}\Big);$$

\noindent
(b) $${{q^{(r-1)(g-1)}}\over{\gamma_{C,r,{\bf
F}_q}(0)}}=O\Big(1\Big).$$

\noindent
(c) For $0\leq d\leq r(g-1)$,
$${{\alpha_{C,r,{\bf F}_q}(d)}\over
{q^{d/2+r+r^2(g-1)}}}=O(1).$$}

\noindent
{\it Proof.}  Following Harder and Narasimhan [HN],  a result of Siegel on 
quadratic forms, which is
equivalent to the fact that Tamagawa number of ${\rm SL}_r$ is 1, may be 
understood via the
following relation on automorphism groups of rank $r$ vector bundles:
$$\sum_{V:r(V)=r,{\rm det}(V)=L}{1\over {\#{\rm 
Aut}(V)}}={{q^{(r^2-1)(g-1)}}\over {q-1}}\cdot
\zeta_C(2)\dots\zeta_C(r).$$ Here $V$ runs over all rank $r$ vector bundles 
with determinant $L$ and
$\zeta_C(s)$ denotes the Artin zeta function of $C$.
Thus, $$0<\beta_{C,r,{\bf F}_q}(L)\leq {{q^{(r^2-1)(g-1)}}\over {q-1}}\cdot
\zeta_C(2)\dots\zeta_C(r).$$ This implies $$\beta_{C,r,{\bf
F}_q}(d)=\prod_{i=1}^{2g}(1-\omega_{C,1,{\bf F}_q}(i))\cdot \beta_{C,r,{\bf 
F}_q}(L)
\leq \prod_{i=1}^{2g}(1-\omega_{C,1,{\bf 
F}_q}(i))\cdot{{q^{(r^2-1)(g-1)}}\over {q-1}}\cdot
\zeta_C(2)\dots\zeta_C(r).$$ Here two facts are used:

\noindent
(1) The number of ${\bf F}_q$-rational points of degree $d$ Jacobian 
$J^d(C)$ is equal to
$\prod_{i=1}^{2g}(1-\omega_{C,1,{\bf F}_q}(i))$; and

\noindent
(2) a result of Desale and Ramanan, which says that for any two $L,L'\in 
{\rm Pic}^d(C)$,
$\beta_{C,r,{\bf F}_q}(L)=\beta_{C,r,{\bf F}_q}(L')$. (See e.g., [DR, Prop 
1.7.(i)])

Thus by Weil's theorem on Riemann Hypothesis on Artin zeta functions 
([W]),  $$|\omega_{C,1,{\bf
F}_q}(i)|=O(q^{1/2}), \qquad i=0,\dots,2g.$$ This then completes the proof 
of (a).

To prove (b), first  note that (b) is equivalent to that, asymptotically, 
the lower bound of
$\gamma_{C,r,{\bf F}_q}(0)$ is at least $q^{(r-1)(g-1)}$. To show this, 
note that
$$\eqalign{\gamma_{C,r,{\bf F}_q}(0)&\geq\sum_{V={\cal O}_C\oplus 
L_2\oplus\dots\oplus L_r,
L_2,\dots,L_r\in {\rm Pic}^0(C), \#\{{\cal 
O}_C,L_2,\dots,L_r\}=r}{{q^{h^0(C,V)}-1}\over {\#{\rm
Aut}(V)}}\cr
  =&{1\over {(q-1)^{r-1}}}\sum_{V={\cal O}_C\oplus L_2\oplus\dots\oplus L_r,
L_2,\dots,L_r\in {\rm Pic}^0(C), \#\{{\cal O}_C,L_2,\dots,L_r\}=r} 1.\cr}$$ 
Now, by
the above cited result of Weil again, as $q\to\infty$, $$\sum_{V={\cal 
O}_C\oplus L_2\oplus\dots\oplus L_r,
L_2,\dots,L_r\in {\rm Pic}^0(C), \#\{{\cal O}_C,L_2,\dots,L_r\}=r} 
1=O(q^{g(r-1)}).$$
So we have (b) as well.

Just as (a), (c) is about to give an upper bound for $\alpha_{C,r,{\bf 
F}_q}(d)$ for $0\leq d\leq
r(2g-2)$. For this, we first recall the following Clifford Lemma.

\noindent
{\bf Clifford Lemma.} (See e.g., [B-PBGN, Theorem
2.1]) {\it Let $V$ be a semi-stable bundle of rank
$r$ and degree
$d$ with $0\leq\mu(V)\leq 2g-2$. Then $$h^0(C,V)\leq r+{d\over 2}.$$}
 
Thus, $$\alpha_{C,r,{\bf F}_q}(d)\leq q^{{d\over 2}+r}\cdot\beta_{C,r,{\bf 
F}_q}(d).$$
With this,  (c) is a direct consequence of (a).

\vskip 0.30cm
\noindent
{\bf 2.1.3. Ugly Formula}
\vskip 0.30cm
Recall that the rationality of $\zeta_{C,r,{\bf F}_q}(s)$ says that there 
exists a degree $2rg$
polynomial $P_{C,r,{\bf F}_q}(t)\in {\bf Q}[t]$ such that
$$Z_{C,r,{\bf F}_q}(t)={{P_{C,r,{\bf F}_q}(t)}\over {(1-t^r)(1-q^rt^r)}}.$$
Thus we may set $$P_{C,r,{\bf F}_q}(t)=\sum_{i=0}^{2rg}a_{C,r,{\bf 
F}_q}(i)t^i.$$ On the other
hand, by the functional equation for $\xi_{C,r,{\bf F}_q}(t)(s)$, we have
$$P_{C,r,{\bf F}_q}(t)=P_{C,r,{\bf F}_q}({1\over {qt}})\cdot q^{rg}\cdot 
t^{2rg}.$$ Thus
by comparing coefficients on both sides, we get the following

\noindent
{\bf Lemma.} {\it With the same notation as above, for $i=0,1,\dots,rg-1$,
$$a_{C,r,{\bf F}_q}(2rg-i)=a_{C,r,{\bf F}_q}(i)\cdot q^{rg-i}.$$}

Now, to determine $P_{C,r,{\bf F}_q}(t)$ and hence $\zeta_{C,r,{\bf 
F}_q}(s)$ it suffices to find
$a_{C,r,{\bf F}_q}(i)$ for $i=0,1,\dots ,rg$.

\noindent
{\bf Proposition.} ({\bf An Ugly Formula}) {\it With the same notation as 
above,
$$\eqalign{~&a_{C,r,{\bf F}_q}(i)\cr
=&\cases{\alpha_{C,r,{\bf
F}_q}(d)-\beta_{C,r,{\bf F}_q}(d),& if\ $0\leq i\leq r-1$;\cr
\alpha_{C,r,{\bf F}_q}(d)
-(q^r+1)\alpha_{C,r,{\bf F}_q}(d-r)+q^r\beta_{C,r,{\bf
F}_q}(d-r),& if\ $r\leq i\leq 2r-1$;\cr
\alpha_{C,r,{\bf F}_q}(d)
-(q^r+1)\alpha_{C,r,{\bf F}_q}(d-r)+q^r\alpha_{C,r,{\bf
F}_q}(d-2r),& if\ $2r\leq i\leq r(g-1)-1$;\cr
-(q^r+1)\alpha_{C,r,{\bf F}_q}(r(g-2))
+q^r\alpha_{C,r,{\bf F}_q}(r(g-3))
+\alpha_{C,r,{\bf F}_q}(r(g-1)),& if\ $i=r(g-1)$;\cr
\alpha_{C,r,{\bf F}_q}(d)
-(q^r+1)\alpha_{C,r,{\bf F}_q}(d-r)+\alpha_{C,r,{\bf
F}_q}(d-2r)q^r,& if\ $r(g-1)+1\leq i\leq rg-1$;\cr
2q^r\alpha_{C,r,{\bf F}_q}(r(g-2))
-(q^r+1)\alpha_{C,r,{\bf F}_q}(r(g-1)),& if\ $i=rg$;\cr}\cr}$$}

\noindent
{\it Proof.} By definition,
$$\eqalign{~&Z_{C,r,{\bf F}_q}(t)\cr
=&(\sum_{d=0}^{r(2g-2)}+\sum_{d=r(2g-2)+1}^\infty)\sum_{V\in [V]\in
{\cal M}_{C,r}(d),d\geq 0} {{q^{h^0(C,V)}-1}\over{\#{\rm Aut}(V)}}t^d\cr
=&\sum_{d=0}^{r(2g-2)}\sum_{V\in [V]\in {\cal
M}_{C,r}(d)} {{q^{h^0(C,V)}-1}\over{\#{\rm Aut}(V)}}t^d\cr
&\qquad+\sum_{i=1}^r\sum_{n=2g-2}^\infty\sum_{d=rn+i} \sum_{V\in [V]\in {\cal
M}_{C,r}(d)} {{q^{rn+i-r(g-1)}-1}\over{\#{\rm Aut}(V)}}t^{rn+i}\cr
=&\sum_{d=0}^{r(2g-2)}\sum_{V\in [V]\in {\cal
M}_{C,r}(d)({\bf F}_q)} {{q^{h^0(C,V)}-1}\over{\#{\rm Aut}(V)}}t^d\cr
&\qquad+{{q^{r(1-g)}}\over{1-(qt)^r}}(qt)^{r(2g-2)}
  \sum_{i=1}^r \beta_{C,r,{\bf F}_q}(i) (qt)^i
-{1\over{1-t^r}}t^{r(2g-2)}\sum_{i=1}^r \beta_{C,r,{\bf F}_q}(i) t^i,\cr}$$
by a similar calculation as in the proof of Lemma 1.2.4.(b).
Now
$$\sum_{d=0}^{r(2g-2)}=\sum_{d=0,r(2g-2)}+\sum_{d=1,r(2g-2)-1}+\dots+\sum_{d 
=r(g-1)-1,r(g-1)+1}+\sum_{d=r(g-1)}.$$
Thus, by Riemann-Roch,  Serre duality and Lemma 2.1.1, we conclude that
$$\eqalign{~&\sum_{d=0}^{r(2g-2)}\sum_{V\in [V]\in {\cal
M}_{C,r}(d)} {{q^{h^0(C,V)}-1}\over{\#{\rm Aut}(V)}}t^d\cr
=&\sum_{d=0}^{r(g-1)-1}\Big[\alpha_{C,r,{\bf F}_q}(d)\Big(
t^d+q^{r(g-1)-d}t^{r(2g-2)-d}\Big)-\beta_{C,r,{\bf
F}_q}(d)\Big(t^d+t^{r(2g-2)-d}\Big)\Big]\cr
&\qquad\qquad+\Big(\alpha_{C,r,{\bf
F}_q}(r(g-1))-\beta_{C,r,{\bf
F}_q}(r(g-1))\Big)\cdot t^{r(g-1)}.\cr}$$ With all this, together with 
Lemmas 1.2.2.2
and 2.1.1, by a couple of pages routine calculation, we are lead to the 
ugly yet very precise
formula in the proposition.
\vskip 0.30cm
\centerline {\bf 2.2. Global Non-Abelian Zeta Functions for Curves}
\vskip 0.30cm
\noindent
{\bf 2.2.1. Definition.} Let ${\cal C}$ be a regular, reduced, irreducible 
projective curve of genus
$g$ defined over a number field $F$. Let $S_{\rm bad}$ be the collection of 
all infinite places and
these finite places of $F$ at which ${\cal C}$  do not have good 
reductions. As usual, a place $v$ of
$F$ is called good if $v\not\in S_{\rm bad}$.
 
Thus, in particular, for any good place $v$ of $F$,   the $v$-reduction of 
${\cal C}$, denoted as
${\cal C}_v$, gives a regular, reduced, irreducible projective curve 
defined over the residue field
field $F(v)$ of $F$ at $v$. Denote the cardinal number of $F(v)$ by $q_v$. 
Then,  by the
construction of Chapter 1, we obtain the associated rank $r$ non-abelian 
zeta function $\zeta_{{\cal
C}_v,r,{\bf F}_{q_v}}(s)$. Moreover, from the rationality of $\zeta_{{\cal
C}_v,r,{\bf F}_{q_v}}(s)$, there exists a degree $2rg$ polynomial $P_{{\cal 
C}_v,r,{\bf F}_{q_v}}(t)\in {\bf Q}[t]$
such that
$$Z_{{\cal C}_v,r,{\bf F}_{q_v}}(t)={{P_{{\cal C}_v,r,{\bf 
F}_{q_v}}(t)}\over {(1-t^r)(1-q^rt^r)}}.$$
Clearly, $$P_{{\cal C}_v,r,{\bf F}_{q_v}}(0)=\gamma_{{\cal C}_v,r,{\bf 
F}_{q_v}}(0)\not=0.$$ Thus it
makes sense to introduce the polynomial $\tilde P_{{\cal C}_v,r,{\bf 
F}_{q_v}}(t)$ with constant
term 1 by setting
$$\tilde P_{{\cal C}_v,r,F(v)}(t):={{P_{{\cal C}_v,r,F(v)}(t)}\over 
{P_{{\cal C}_v,r,F(v)}(0)}}.$$
Now by definition, {\it the rank $r$ non-abelian
zeta function $\zeta_{{\cal C},r,F}(s)$ of} ${\cal C}$ over $F$ is the 
following Euler product
$$\zeta_{{\cal C},r,F}(s)
:=\prod_{v:{\rm good}}{1\over{
\tilde P_{{\cal C}_v,r,{\bf F}_{q_v}}(q_v^{-s})}},\hskip 2.0cm {\rm 
Re}(s)>>0.$$

Clearly, when $r=1$, $\zeta_{{\cal C},r,F}(s)$ coincides with the classical 
Hasse-Weil zeta function
for $C$ over $F$ ([H]).

\noindent
{\bf 2.2.2. Convergence.} At this earlier stage of the study of our 
non-abelian zeta functions, the
central problem is to justify the above definition. That is to say, to show 
indeed the Euler product
in 2.2.1 converges. In this direction, we have the following
\vskip 0.30cm
\noindent
{\bf Conjecture.} {\it Let ${\cal C}$ be a regular, reduced, irreducible 
projective curve of genus
$g$ defined over a number field $F$. Then its associated rank $r$ global 
non-abelian zeta function
$\zeta_{{\cal C},r,F}(s)$  admits a meromorphic continuation to the whole 
complex $s$-plane.}

Recall that even  when $r=1$, i.e., for the classical Hasse-Weil zeta 
functions, this conjecture has
not been confirmed. However, for general $r$, we have the following
\vskip 0.30cm
\noindent
{\bf Theorem.} {\it  Let ${\cal C}$ be a regular, reduced, irreducible 
projective curve defined over
a number  field $F$. Then its associated rank r global non-abelian zeta 
function
$\zeta_{{\cal C},r,F}(s)$ converges  when ${\rm Re}(s)\geq 1+g+(r^2-r)(g-1)$.}
\vskip 0.30cm
\noindent
{\it Proof.} Clearly, it suffices to show that for the reciprocal roots 
$\omega_{C,r,{\bf F}_q}(i),
i=1,\dots,2rg$ of $P_{C,r,{\bf F}_q}(t)$ associated to curves $C$ over 
finite fields ${\bf F}_q$,
$$|\omega_{C,r,{\bf F}_q}(i)|=O(q^{g+(r^2-1)(g-1)}).$$
Thus we are lead to estimate coefficients of $P_{C,r,{\bf F}_q}(t)$. Since 
we have the ugly yet very
precise formula for these coefficients, i.e., Lemma and Proposition 2.1.3, 
it suffices to give upper
bounds for $\alpha_{C,f,{\bf F}_q}(i), \beta_{C,f,{\bf F}_q}(j)$ and a 
lower bound
for $\gamma_{C,f,{\bf F}_q}(0)$, the constant term of $P_{C,r,{\cal 
F}_q}(t)$. This is achieved in
Proposition 2.1.2.

\noindent
{\bf 2.2.3. Working Hypothesis.} Like in the theory for abelian zeta 
functions, we want to use our
non-abelian zeta functions  to study  non-abelian aspect of
arithmetic of curves. Motivated by the classical analytic class number 
formula for Dedekind zeta
functions and its counterpart BSD conjecture for Hasse-Weil zeta functions 
of elliptic curves, we
expect that our non-abelian zeta function could be used to understand the 
Weil-Petersson volumes of
moduli space of stable bundles as well as the associated Tamagawa measures.

For doing so, we then also need to introduce the local factors for \lq bad' 
places. This may be done
as follows: for $\Gamma$-factors, we take these coming from the functional 
equation for
  $\zeta_F(rs)\cdot\zeta_F(r(s-1))$, where $\zeta_F(s)$ denotes
the standard Dedekind zeta function for $F$; while for finite bad places, 
we do as follows: first,
use the semi-stable reduction for curves to find a semi-stable model for 
${\cal C}$, then use
Seshadri's moduli spaces of parabolic bundles to construct polynomials for 
singular fibers, which
usually have degree lower than $2rg$. With all this being done, we then can 
introduce the so-called
completed rank
$r$ non-abelian zeta function  for
${\cal C}$ over
$F$, or better,  the  completed rank $r$ non-abelian zeta function 
$\xi_{X,r,{\cal O}_F}(s)$ for a
semi-stable model $X\to {\rm Spec}({\cal O}_F)$ of ${\cal C}$. Here ${\cal 
O}_F$ denotes the ring of
integers of $F$. (If necessary, we take a finite extension of $F$.)

\noindent
{\bf Conjecture.} {\it $\xi_{X,r,{\cal O}_F}(s)$ is holomorphic and 
satisfies the functional equation
$$\xi_{X,r,{\cal O}_F}(s)=\pm\, \xi_{X,r,{\cal O}_F}(1+{1\over r}-s).$$}

Moreover, we expect that for certain classes of curves, the inverse Mellin 
transform of our
non-abelian zeta functions are naturally associated to certain modular 
forms of weight $1+{1\over
r}$.

\noindent
{\it Remark.} From our study for non-abelian zeta functions of elliptic 
curves [We2], we obtain the
following so-called \lq absolute Euler product' for rank 2 zeta functions 
of elliptic curves
$$\eqalign{\zeta_2(s)=&\prod_{p>2; {\rm prime}}{1\over
{1+(p-1)p^{-s}+(2p-4)p^{-2s}+(p^2-p)p^{-3s}+p^2p^{-4s}}}\cr
=&\prod_{p>2; {\rm prime}}{1\over {A_p(s)+B_p(s)p^{-2s}}},\qquad {\rm 
Re}(s)>2\cr}$$
with $$A_p(s)=1+(p-1)p^{-s}+(p-2)p^{-2s},\qquad 
B_p(s)=(p-2)+(p^2-p)p^{-s}+p^2p^{-2s}.$$
Set $t:=q^{-s}$ and $a_p(t):=A_p(s), b_p(t):=B_p(s)$. Then in ${\bf Z}[t]$, 
we have the factorization
$$a_p(t)=(1+(p-2)t)(1+t),\qquad
b_p(t)=((p-2)+pt)(1+pt)$$ and
$$a_p({1\over {pt}})={1\over {p^2t^2}}\cdot b_p(t).$$
We hope a similar discussion on Hecke operators as in the classical theory 
of $L$-functions, or
better modular forms, works here as well. In particular we ask the following

\noindent
{\bf Questions.} {\it (a) What can we say about $\zeta_2(s)$?

\noindent
(b) What is the importance of the factorization of $a_p(t)$ and $b_p(t)$ in 
${\bf Z}[t]$?}
\vskip 0.45cm
\centerline {\li Chapter 3. Non-Abelian Zeta
Functions and Infinitesimal Structures of Brill-Noether Loci}
\vskip 0.30cm
In this chapter, we study the infinitesimal structures of the so-called 
non-abelian Brill-Noether
loci for rank two semi-stable vector bundles over genus two curves. As an 
application, we calculate
the corresponding rank two non-abelian zeta functions for genus two curves. 
During this process, we
see clearly how Weierstrass points, intrinsic arithmetic invariants of 
curves [We3], contribute to
our zeta functions among others.

We in this chapter assume that the characteristic of the base field is 
strictly bigger than 2 for
simplicity.
\vskip 0.30cm
\centerline {\bf 3.1. Infinitesimal Structures of Non-Abelian Brill-Noether 
Loci}
\vskip 0.30cm
\noindent
{\bf 3.1.1. Invariants $\beta_{C,2, {\bf F}_q}(d)$.} Let $C$ be a regular 
reduced irreducible
projective curve defined over ${\bf F}_q$. Here we want to calculate 
Harder-Narasimhan numbers
$\beta_{C,2,{\bf F}_q}(d)$ for all $d$. Note that from Lemma 2.1.1, 
$$\beta_{C,2,{\bf F}_q}(d)
=\beta_{C,2,{\bf F}_q}(d+2n).$$ So it suffices to calculate 
$\beta_{C,2,{\bf F}_q}(d)$ when $d=0,1$.
For this, we cite the following result of Desale and Ramanan:

\noindent
{\bf Proposition.} ([DR]) {\it With the same notation as above, for $L\in 
{\rm Pic}^d(C), d=0,1$,
$$\beta_{C,2,{\bf F}_q}(L)={{q^3}\over
{q-1}}\cdot\zeta_C(2)-q\prod_{i=1}^4(1-\omega_i)\cdot \sum_{d_1+d_2=d, 
d_1>d_2}{{\beta_{C,1,{\bf
F}_q}(d_1)\beta_{C,1,{\bf F}_q}(d_2)}\over {q^{d_1-d_2}}}.$$ Here 
$\zeta_C(s)$ denotes the Artin
zeta function for $C$ and $\omega_1,\dots,\omega_4$ are the roots of the 
associated
$Z$-function $Z_C(s)$, i.e., $\omega_{C,1,{\bf F}_q}(i), i=0,\dots, 
4=2\times 2$ in our notation.}

Thus, in particular, $\beta_{C,2,{\bf F}_q}(L)$ is independent of $L$.

\noindent
{\bf Lemma.} {\it With the same notation as above, for $d=0,1$
$$\beta_{C,2,{\bf F}_q}(d)={{q^3}\over {q-1}}\cdot\zeta_C(2)\cdot 
\prod_{i=1}^4(1-\omega_i)
-{{q^{d+1}}\over {(q-1)^2(q^2-1)}}\cdot \prod_{i=1}^4(1-\omega_i)^4.$$}

\noindent
{\it Proof.} This comes from the following two facts:

\noindent
(1) for all $d$, $$\beta_{C,1,{\bf F}_q}(d)={{\prod_{i=1}^4(1-\omega_i)}\over
{q-1}};$$

\noindent
(2) the number of ${\bf F}_q$-rational points of ${\rm Pic}^d(C)$ is equal to
$\prod_{i=1}^4(1-\omega_i)$.

\noindent
{\bf 3.1.2. Infinitesimal structures: a taste.} Here we want to calculate 
$\alpha_{C,2,{\bf
F}_q}(0)$. By Lemma 3.1.1, it suffices to give $\gamma_{C,2,{\bf F}_q}(0)$. 
So we are lead to
study
$\gamma_{C,2,{\bf F}_q}(L)$ which is
supported over the Brill-Noether locus $$W_{C,2}^0(L):=\{[V]\in {\cal 
M}_{C,2}(L):h^0(C,{\rm
Gr}(V))\geq 1\}.$$ (In general, as in the beautiful paper [B-PGN],
we define the Brill-Noether locus
$$W_{C,2}^k(L):=\{[V]\in {\cal M}_{C,2}(L):h^0(C,{\rm Gr}(V))\geq k+1\}.)$$

Note that no degree zero stable bundle admits non-trivial
global sections, so
$W_{C,2}^0(L):=\{[{\cal O}_C\oplus L]\}$ consists of only a single point.

\noindent
(a) If $L={\cal O}_C$, then $W_{C,2}^0({\cal O}_C)=W_{C,2}^1({\cal O}_C).$ 
Moreover, infinitesimally,
$V={\cal O}_C\oplus {\cal O}_C$ or $V$ corresponds to all non-trivial 
extensions
$$0\to {\cal O}_C\to V\to {\cal O}_C\to 0$$ which are parametrized by ${\bf 
P}{\rm Ext}^1({\cal
O}_C,{\cal O}_C)\simeq {\bf P}^1$. Thus, by definition,
$$\gamma_{C,2,{\bf F}_q}({\cal O}_C)={{q^2-1}\over 
{(q^2-1)(q^2-q)}}+(q+1)\cdot {{q-1}\over {q(q-1)}}
={q\over {q-1}}.$$

\noindent
(b) If $L\not ={\cal O}_C$, then,
infinitesimally,
$V={\cal O}_C\oplus L$ or $V$ corresponds to the single non-trivial extension
$$0\to {\cal O}_C\to V\to L\to 0.$$  Thus, by definition,
$$\gamma_{C,2,{\bf F}_q}(L)={{q-1}\over {(q-1)^2}}+{{q-1}\over {q-1}}
={q\over {q-1}}.$$

Thus we have the following

\noindent
{\bf Lemma.} {\it With the same notation as above, for all $L\in {\rm 
Pic}^0(C)$,
$$\gamma_{C,2,{\bf F}_q}(L)={q\over {q-1}}.$$  In particular,
$$\gamma_{C,2,{\bf F}_q}(0)={q\over {q-1}}\cdot \prod_{i=1}^4(1-\omega_i).$$}

\noindent
{\bf 3.1.3. Invariants $\alpha_{C,2{\bf F}_q}(1)$.} As before, it suffices 
to calculate
$\gamma_{C,2{\bf F}_q}(L)$ for all $L\in {\rm Pic}^1(C)$. Note that in this 
case, all bundles are
stable, so ${\rm Aut}(V)\simeq {\bf F}_q^*$ and
$$W_{C,2}^0(L)\simeq\{V:{\rm stable}, r(V)=2,{\rm det}(V)=L, h^0(C,V)\geq
1\}.$$ Moreover, by [B-PGN, Prop. 3.1], $$W_{C,2}^0(L)=\{V:{\rm stable}, 
r(V)=2,{\rm
det}(V)=L,h^0(C,V)=1\}$$ and any $V\in W_{C,2}^0(L)$ admits a non-trivial 
extension
$$0\to {\cal O}_C\to V\to L\to 0.$$ On the other hand, any non-trivial 
extension
$$0\to {\cal O}_C\to V\to L\to 0$$ gives rise a stable bundle. So in fact
$$W_{C,2}^0(L)\simeq {\bf P}{\rm Ext}^1(L,{\cal O}_C)\simeq {\bf P}^1.$$
Thus we have the following

\noindent
{\bf Lemma.} {\it With the same notation as above, for $L\in {\rm Pic}^1(C)$,
$$W_{C,2}^0(L)\simeq {\bf P}^1,\qquad {\rm and}\qquad \gamma_{C,2,{\bf 
F}_q}(L)=q+1.$$
In particular, $$\gamma_{C,2{\bf F}_q}(1)=(q+1)\cdot 
\prod_{i=1}^4(1-\omega_i).$$}

\noindent
{\bf 3.1.4. Non-abelian Brill-Noether loci and their infinitesimal 
structures.} We here want to
calculate $\gamma_{C,r,{\bf F}_q}(2)$. This is the most complicated level 
as $2=r(g-1)$.
For this purpose, we need to understand the structures of the non-abelian 
Brill-Noether loci
$W_{C,2}^0(L)$ and $W_{C,2}^1(L)$ for $L\in {\rm Pic}^2(C)$.

We begin with recalling the structure of the map $\pi:C\times C/S_2\to {\rm 
Pic}^2(C)$. Here $S_2$
denotes the symmetric group of two symbols which acts naturally on $C\times 
C$ via $(x,y)\mapsto
(y,x)$. One checks that $\pi$ is a one point blowing-up centered at the 
canonical line
bundle $K_C$ of $C$. For later use, denote by $\Delta$ the image of the 
diagonal of $C\times C$ in
${\rm Pic}^2(C)$.

Next, we want to understand the structure of sublocus $W_{C,2}^0(L)^{\rm 
ss}$ of $W_{C,2}^0(L)$
consisting of non-stable but semi-stable vector bundles.

By definition, for any $V\in [V]\in  W_{C,2}^0(L)^{\rm ss}$, ${\rm 
Gr}(V)={\cal O}_C(P)\oplus L(-P)$
for a suitable (${\bf F}_q$-rational) point $P\in C$. Thus accordingly,

\noindent
(a) if $L\not=K_C$, then $W_{C,2}^0(L)^{\rm ss}$ is parametrized by (${\bf 
F}_q$-rational points
of) $C$, due to the fact that now $h^0(C,L)=1$. Write also $L={\cal 
O}_C(A+B)$ with two points $A,B$
of $C$, which are unique from the above discussion on the map $\pi$, we 
then conclude that
$$W_{C,2}^1(L)=\{[{\cal O}_C(A)\oplus {\cal O}_C(B)]\}.$$

\noindent
(b) if $L=K_C$, then for any $P$, $K_C={\cal O}_C(P+\iota(P))$ where 
$\iota: C\to C$ denotes the
canonical involution on $C$. So $$W_{C,2}^0(L)^{\rm ss}=\{[{\cal 
O}_C(P)\oplus {\cal
O}_C(\iota(P))]:P\in C\}.$$ Thus, indeed $W_{C,2}^0(L)^{\rm ss}$ is 
parametrized by ${\bf
P}^1$. Moreover, $$W_{C,2}^1(K_C)=W_{C,2}^0(L)^{\rm ss}=\{[{\cal 
O}_C(P)\oplus {\cal
O}_C(\iota(P))]:P\in C\}.$$

On the other hand, it is easily to check that every non-trivial extension 
$$0\to {\cal O}_V\to W\to
L\to 0$$ gives rise to a semi-stable vector bundle $W$, and if $W$ is not 
stable, then there exists
a point $Q\in C$ such that $W$ may also be given by the non-trivial 
extension $$0\to {\cal O}_C(Q)\to
W\to L(-Q)\to 0.$$ Note also that the kernel of the natural map $H^1(C,{\rm 
Hom}(L,{\cal
O}_C))\to H^1(C,{\rm Hom}(L(-Q),{\cal O}_C))$ is one dimensional. So among all
non-trivial extensions $0\to {\cal O}_C\to V\to L\to 0$, which are 
parametrized by ${\bf P}{\rm
Ext}^1(L,{\cal O}_C)\simeq {\bf P}^2$, the non-stable (yet semi-stable) 
vector bundles are
paremetrized by (${\bf F}_q$-rational points of) $C$ when $L\not=K_C$ by 
(a) and ${\bf P}^1$ when
$L=K_C$ by (b) above respectively. (See [NR, Lemma 3.1]) In this way, we 
have proved the
following  result on non-abelian Brill-Noether loci for moduli space of
${\cal M}_{C,2}(L)$ with $L$ a degree 2 line bundles on a genus two curve, 
which is not covered by
[B-PGN]:

\noindent
{\bf Lemma.} {\it With the same notation as above, $W_{C,2}^0(L)\simeq {\bf 
P}{\rm Ext}^1(L,{\cal
O}_C)\simeq {\bf P}^2,$ in which the locus $W_{C,2}^0(L)^{\rm ss}$ of 
semi-stable but not stable
bundles is parametrized by
$C$ and ${\bf P}^1$  in the cases  $L\not= K_C$ and $L=K_C$ respectively. 
More precisely,

\noindent
(a) if $L={\cal O}_C(A+B)\not=K_C$ with $A,B$ two points of $C$, then 
$W_{C,2}^0(L)^{\rm ss}$, as a
birational image of $C$ under the complete linear system $K_C(A+B)$, is a 
degree 4 plane curve
with a single node located at $W_{C,2}^1(L)=\{[{\cal O}_C(A)\oplus {\cal 
O}_C(B)]\}$;

\noindent
(b) If $L=K_C$, as a degree 2 regular plane curve,
$$W_{C,2}^1(K_C)=W_{C,2}^0(L)^{\rm ss}=\{[{\cal O}_C(P)\oplus {\cal 
O}_C(\iota(P))]:P\in C\}\simeq
{\bf P}^1.$$}

Next, we study the infinitesimal structures of non-abelian Brill-Noether 
loci. Set
$$W_{C,2}^0(L)^{s}:=W_{C,2}^0(L)\backslash W_{C,2}^0(L)^{\rm ss}.$$
Then the infinitesimal structure of $W_{C,2}^0(L)$ at points $[V]\in 
W_{C,2}^0(L)^{s}$ is simple:
each $[V]$ consists a single stable rank two vector bundle with ${\rm 
det}(V)=L, h^0(C,V)=1$ and
${\rm Aut}(V)\simeq {\bf F}_q^*$.

Now we consider $W_{C,2}^0(L)^{\rm ss}$.

\noindent
(a)  $L\not=K_C$. Then there exists two points $A,B$ of $C$ such that 
$L={\cal O}_C(A+B)$.
Thus, for any $V\in [{\cal O}_C(P)\oplus {\cal O}_C(A+B-P)]\not\in 
W_{C,2}^1(L)$,
$V$ is given by an extension $0\to {\cal O}_C(P)\to V\to {\cal 
O}_C(A+B-P)\to 0$ due to the fact that
for the non-trivial extension $0\to {\cal O}_C(A+B-P)\to W\to {\cal 
O}_C(P)\to 0$, $h^0(C,W)=0$.
Thus,  each class $[{\cal O}_C(P)\oplus {\cal O}_C(A+B-P)]\not\in W_{C,2}^1(L)$
consists of exact two vector bundles, i.e.,
$V_1={\cal O}_C(P)\oplus {\cal O}_C(A+B-P)$ and $V_2$ given by the 
non-trivial extension
$0\to {\cal O}_C(P)\to V\to {\cal O}_C(A+B-P)\to 0$. Clearly, 
$h^0(C,V_1)=h^0(C,V_2)=1$ and $\#{\rm
Aut}(V_1)=(q-1)^2,  \#{\rm Aut}(V_2)=q-1;$

  To study $W_{C,2}^1(L)=\{[{\cal O}_C(A)\oplus {\cal O}_C(B)]\}$, we divide 
it into  two
subcases.

\noindent
(i)  $A\not=B$. Then there are exact three vector bundles in the class 
$[{\cal O}_C(A)\oplus {\cal
O}_C(B)]$. They are $V_0={\cal O}_C(A)\oplus {\cal O}_C(B)$,
$V_1$ given by the non-trivial extension $0\to {\cal O}_C(A)\to V_2\to 
{\cal O}_C(B)\to 0$ and $V_2$
given by the non-trivial extension $0\to {\cal O}_C(B)\to V_2\to {\cal 
O}_C(A)\to 0$. Clearly,
$h^0(C,V_1)=2, h^0(C,V_1)=h^0(C,V_2)=1$ and $\#{\rm Aut}(V_0)=(q-1)^2,\
\#{\rm Aut}(V_1)=\#{\rm Aut}(V_2)=q-1$;

Thus in particular, $$\gamma_{C,2,{\bf F}_q}(L)=(q^2+q+1-(N_1-1))\cdot 
{{q-1}\over {q-1}}
+(N_1-2)\Big({{q-1}\over {(q-1)^2}}+{{q-1}\over 
{q-1}}\Big)+\Big({{q^2-1}\over {(q-1)^2}}+
{{q-1}\over {q-1}}+{{q-1}\over {q-1}}\Big).$$ Here 
$N_1=q+1-(\omega_1+\dots+\omega_4)$ denotes the
number of ${\bf F}_q$-rational points of $C$.

\noindent
(ii) $A=B$. Then the infinitesimal structure at $[{\cal O}_C(A)\oplus {\cal
O}_C(A)]$ is as follows: an independent point corresponding to $V_0={\cal 
O}_C(A)\oplus {\cal
O}_C(A)$ and a projective line parametrizing all non-trivial extension 
$0\to {\cal O}_C(A)\to V\to
{\cal O}_C(A)\to 0$. Clearly,
$h^0(C,V_0)=2, h^0(C,V)=1$ and $\#{\rm Aut}(V_0)=(q^2-1)(q^2-q),\ \#{\rm 
Aut}(V)=q(q-1);$

Thus in particular, $$\gamma_{C,2,{\bf F}_q}(L)=(q^2+q+1-(N_1-1))\cdot 
{{q-1}\over {q-1}}
+(N_1-2)\Big({{q-1}\over {(q-1)^2}}+{{q-1}\over {q-1}}\Big)+\Big({{q^2-1}\over
{(q^2-1)(q^2-q)}}+(q+1){{q-1}\over {q(q-1)}}\Big).$$

\noindent
(b) $L=K_C$. Then $K_C={\cal O}_C(P+\iota(P))$ for all points $P$. 
Therefore for all $[V]\in
W_{C,2}^1(L)=W_{C,2}^0(L)^{\rm ss}$, $[V]=[{\cal O}_C(P)\oplus {\cal 
O}_C(\iota(P))]$. Accordingly,
two subcases:

\noindent
(i) $P\not=\iota P$. Then there are exact three vector bundles in the class 
$[{\cal O}_C(P)\oplus
{\cal O}_C(\iota(P)]$. They are $V_0={\cal O}_C(P)\oplus {\cal O}_C(\iota(P))$,
$V_1$ given by the non-trivial extension $0\to {\cal O}_C(P)\to V_1\to 
{\cal O}_C(\iota(P))\to 0$
and $V_2$  given by the non-trivial extension $0\to {\cal O}_C(\iota(P))\to 
V_2\to {\cal O}_C(P)\to
0$. Clearly,
$h^0(C,V_0)=2, h^0(C,V_1)=h^0(C,V_2)=1$ and $\#{\rm Aut}(V_0)=(q-1)^2,\
\#{\rm Aut}(V_1)=\#{\rm Aut}(V_2)=q-1$;

\noindent
(ii) $P=\iota(P)$ a Weierstrass point, all of which are four. Then the 
infinitesimal structure at
$[{\cal O}_C(P)\oplus {\cal O}_C(P)]$ is as follows: an independent point 
corresponding to $V_0={\cal
O}_C(P)\oplus {\cal O}_C(P)$ and a projective line parametrizing all 
non-trivial extension $0\to
{\cal O}_C(P)\to V\to {\cal O}_C(P)\to 0$. Clearly,
$h^0(C,V_0)=2, h^0(C,V)=1$ and $\#{\rm Aut}(V_0)=(q^2-1)(q^2-q),\ \#{\rm 
Aut}(V)=q(q-1).$

Thus, in particular,
$$\eqalign{~&\gamma_{C,2,{\bf F}_q}(K_C)=(q^2+q+1-(q+1))\cdot {{q-1}\over
{q-1}}\cr
&\qquad+(q+1-4)\Big({{q^2-1}\over {(q-1)^2}}+
{{q-1}\over {q-1}}+{{q-1}\over {q-1}}\Big)+4\Big({{q^2-1}\over 
{(q^2-1)(q^2-q)}}+(q+1){{q-1}\over
{q(q-1)}}\Big).\cr}$$

So we have completed the proof of the following

\noindent
{\bf Proposition.} {\it With the same notation as above,

\noindent
(a) For $L\not=K_C$,

\noindent
(i) if $L\not\in \Delta$, $\gamma_{C,2,{\bf F}_q}(L)={{q^3+2q-3+N_1}\over 
{q-1}}$;

\noindent
(ii) if $L\in \Delta$, $\gamma_{C,2,{\bf F}_q}(L)={{q^3-2+N_1}\over {q-1}}$;

\noindent
(b) if $L=K_C$, $\gamma_{C,2,{\bf F}_q}(L)=q^2+3q-3$.

\noindent
In particular, $$\gamma_{C,2,{\bf 
F}_q}(2)=\Big(\prod_{i=1}^4(1-\omega_i)-(q+1)\Big)\cdot
{{q^3+2q-3+N_1}\over {q-1}}+q\cdot {{q^3-2+N_1}\over {q-1}}+q^2+3q-3.$$}

In this way, by using the ugly formula in Prop 2.1.3, we can finally write 
down the rank two
non-abelian zeta functions for genus two curves, where a degree 8 
polynomial is involved. We leave
this to the reader.
\vskip 0.45cm
\centerline {\bf REFERENCES}
\vskip 0.30cm
\item{[A]} E. Artin, Quadratische K\"orper im Gebiete der h\"oheren
Kongruenzen, I,II, {\it Math. Zeit}, {\bf 19} 153-246 (1924) (See also
{\it Collected Papers}, pp. 1-94,  Addison-Wesley 1965)
\vskip 0.30cm
\item{[B-PGN]} L. Brambila-Paz, I. Grzegorczyk, \& P.E. Newstead, Geography 
of Brill-Noether Loci
for Small slops, J. Alg. Geo. {\bf 6}  645-669 (1997)
\vskip 0.30cm
\item{[DR]} U.V. Desale \& S. Ramanan, Poincar\'e polynomials of the variety of
stable bundles, Maeh. Ann {\bf 216}, 233-244 (1975)
\vskip 0.30cm
\item{[HN]} G. Harder \& M.S. Narasimhan, On the cohomology groups of 
moduli spaces
of vector bundles over curves, Math Ann. {\bf 212}, (1975) 215-248
\vskip 0.30cm
\item{[H]} H. Hasse, {\it Mathematische Abhandlungen}, Walter
de Gruyter, Berlin-New York, 1975.
\vskip 0.30cm
\item{[M]} D. Mumford, {\it Geometric Invariant Theory}, Springer-Verlag, 
Berlin
(1965)
\vskip 0.30cm
\item{[NR]} M.S. Narasimhan \& S. Ramanan, Moduli of vector bundles on a 
compact Riemann surfaces,
Ann. of Math. {\bf 89} 14-51 (1969)
\vskip 0.30cm
\item {[Se]} C. S. Seshadri, {\it Fibr\'es vectoriels sur les courbes 
alg\'ebriques}, Asterisque {\bf
96}, 1982
\vskip 0.30cm
\item{[W]} A. Weil, {\it Sur les courbes alg\'ebriques et les vari\'et\'es qui
s'en d\'eduisent}, Herman, Paris (1948)
\vskip 0.30cm
\item{[We1]} L. Weng, Riemann-Roch Theorem, Stability, and New Zeta
Functions for Number Fields, preprint, math.AG/0007146
\vskip 0.30cm
\item{[We2]} L. Weng,  Refined Brill-Noether Locus and Non-Abelian Zeta
Functions  for Elliptic  Curves, preprint, math.AG/0101183
\vskip 0.30cm
\item{[We3]} L. Weng, Weierstrass Groups, an appendix to [We2]\end